\documentclass[12pt]{amsart}
\usepackage{amscd}     
\usepackage{amssymb}
\usepackage{amsmath, amsthm, graphics}
\usepackage{xypic}     
\LaTeXdiagrams        
\usepackage[all]{xy}
\xyoption{2cell} \UseAllTwocells \xyoption{frame} \CompileMatrices
\allowdisplaybreaks[3]
\usepackage{amsfonts}
\usepackage[normalem]{ulem}
\usepackage{hyperref}
\usepackage{tikz}
\usepackage{mathtools}
\usepackage{verbatim} 

\usepackage[margin=4cm]{geometry}

\usepackage{latexsym}
\usepackage{epsfig}
\usepackage{amsfonts}
\usepackage{enumerate}
\usepackage{times,float}
\usepackage{graphics}

\usepackage[dvipsnames]{xcolor}

\usepackage{setspace}

\newtheorem{theorem}{Theorem}[section]

\newtheorem{definition}[theorem]{Definition}

\theoremstyle{remark}

\theoremstyle{remark}

\newtheorem{remark}[subsection]{Remark}

\numberwithin{equation}{section}

\title{Descendant and Fourier-Mukai equivalences for simple flops}
\author{Jiun-Cheng Chen}
\address{Department of Mathematics\\ Third General Building\\ National Tsing-Hua University\\ No. 101 Sec. 2 Kuang Fu Road\\ Hsinchu, 30043\\ Taiwan}
\email{jcchen@math.nthu.edu.tw}

\author[Hsian-Hua Tseng]{Hsian-Hua Tseng}
\address{Department of Mathematics\\ Ohio State University\\ 100 Math Tower, 231 West 18th Ave.\\ Columbus, OH 43210\\ USA}
\email{hhtseng@math.ohio-state.edu}
\date{\today}

\begin{document}

\onehalfspacing

\begin{abstract}
For a simple flop $X\dashrightarrow X'$, we construct a correspondence between genus $0$ descendant Gromov-Witten theories of $X$ and $X'$. We show that the Fourier-Mukai equivalence induced by $X\dashrightarrow X'$ is compatible, in a precise sense, with the descendant correspondence.
\end{abstract}

\maketitle

{\centering \em Dedicated to Professor Jun Li on the occasion of his $1000001^{th}$ birthday\par}

\setcounter{tocdepth}{1}
\tableofcontents

\setcounter{section}{-1}

\section{Introduction}
Consider a simple flop $$X\dashrightarrow X'.$$ 
Here $X$ and $X'$ are smooth projective varieties over $\mathbb{C}$, and there is a positive integer $r$ such that: 
\begin{enumerate}
    \item 
The exceptional locus of $X\dashrightarrow X'$ is $\mathbb{P}^r\simeq Z\subset X$, with normal bundle $N_{Z/X}\simeq \mathcal{O}_{\mathbb{P}^r}(-1)^{\oplus r+1}$.
\item 
The blowup $\tilde{X}=\text{Bl}_ZX\xrightarrow{\pi_X} X$ has exceptional divisor $\mathbb{P}^r\times \mathbb{P}^r\simeq \tilde{Z}=\mathbb{P}(N_{Z/X})\to Z\simeq \mathbb{P}^r$.
\item 
The contraction of the other $\mathbb{P}^r$ in $E$ gives $\tilde{X}\xrightarrow{\pi_{X'}} X'$, where $E$ contracts to $\mathbb{P}^r\simeq Z'\subset X'$ with normal bundle $N_{Z'/X'}\simeq \mathcal{O}_{\mathbb{P}^r}(-1)^{\oplus r+1}$. Also, $\tilde{X}=\text{Bl}_{Z'}X'$ and $\pi_{X'}$ is the blowup map $\text{Bl}_{Z'}X'\to X'$. 
\end{enumerate}
The construction may be summarized in the following diagram:
\begin{equation}
\xymatrix{
\,&\tilde{Z}\ar@{^{(}->}[r]\ar[dl]\ar[dr]|\hole & \tilde{X}\ar[dl]|-(0.3){\pi_X}\ar[dr]^{\pi_{X'}} & \,\\
Z\ar@{^{(}->}[r]\ar[dr]&X\ar[dr]\ar@[violet]@{-->}@/_0.49pc/[rr] & Z'\ar@{^{(}->}[r]\ar[dl]|\hole & X'\ar[dl]\\
\,&\text{pt}\ar@{^{(}->}[r] & \bar{X}. & \,
}    
\end{equation}
Details of the construction of simple flops can be found in \cite{LLW}.

The goal of this paper is to construct the following diagram and prove its commutativity:
\begin{equation}\label{diagram:descendant_FM}
\xymatrix{
K(X)\ar[r]^{\mathbb{FM}}\ar[d]^{\Psi^X} & K(X')\ar[d]^{\Psi^{X'}}\\
{\widetilde{\mathcal{H}}^{X}}\ar[r]^{\mathbb{U}} & {\widetilde{\mathcal{H}}^{X'}}.
}    
\end{equation}

The top arrow $\mathbb{FM}$ is the homomorphism between $K$-groups\footnote{In our setting, the $K$-group of locally free sheaves is isomorphic to the Grothendieck group of coherent sheaves, because resolution property holds for smooth projective varieties.} induced from the Fourier-Mukai equivalence 
$$R(p_{X'})_*Lp_X^*: D^b(X)\overset{\simeq}{\longrightarrow} D^b(X'),$$
as proven in \cite{BO}

The vertical arrows $\Psi^X$ and $\Psi^{X'}$ are maps coming from Iritani's integral structure \cite{I}. The definition of these maps can be found in e.g. \cite[Theorem 6.1]{CIJ}. We present the definition for $X$, the definition for $X'$ is parallel.
\begin{definition}
 Define the operators $\text{\em deg}^{X}_0$, $\rho^X$, and $\mu^X$ on $H^*(X)$ as follows. For $\phi\in H^k(X)$, 
\begin{equation*}
\begin{split}
    &\text{\em deg}^{X}_0(\phi)= k\phi\, ,\\
    &\mu^X(\phi)=\left(\frac{k}{2}-\frac{\text{\em dim}\,X}{2}\right)\phi\, ,\\
    &\rho^X(\phi)= c_1(X)\cup\phi\, .
    \end{split}    
\end{equation*} 
The {\em multi-valued Givental space} 
$\widetilde{\mathcal{H}}^X$
for $X$ is defined\footnote{In this paper, we restrict to even cohomology.} by 
$$\widetilde{\mathcal{H}}^X=H^*(X)[[\log(z)]]((z^{-1})).$$
See e.g. \cite{G_symp, G} for more discussions about Givental spaces.

Let $\Psi^X: K(X)\to \widetilde{\mathcal{H}}^X$ be defined by 
\begin{equation*}
  \Psi^X(E)= z^{-\mu^X}z^{\rho^X}\left(\Gamma_X\cup (2\pi \sqrt{-1})^{\frac{\text{\em deg}^{X}_0}{2}}\text{\em ch}(E)\right)\, ,  
\end{equation*}
where $\text{\em ch}(-)$ is the Chern character, $\Gamma_X=\Gamma(T_X)\in H^*(X)$ is the Gamma class of $X$ of \cite{I} and \cite[Section 3.1]{CIJ}, and 
the operators 
$$z^{-\mu^X}:\widetilde{\mathcal{H}}^X\to \widetilde{\mathcal{H}}^X\,,
\ \ \ \ \
z^{\rho^X}:\widetilde{\mathcal{H}}^X\to \widetilde{\mathcal{H}}^X$$ are defined by 
$$z^{-\mu^X}=\sum_{k\geq 0}\frac{\left(-\mu^X\log z\right)^k}{k!}\, ,
\ \ \ \ \
z^{\rho^X}=\sum_{k\geq 0}\frac{\left(\rho^X\log z\right)^k}{k!}\, .$$
\end{definition}

The bottom arrow $\mathbb{U}$ is constructed to equate genus $0$ {\em descendant} Gromov-Witten theories of $X$ and $X'$, based on the equality of primary/ancestor Gromov-Witten theories proven for genus $0$ in \cite{LLW}. See Section \ref{sec:descendant} below.

The main result of this paper is the following
\begin{theorem}\label{thm:comm}
The diagram (\ref{diagram:descendant_FM}) is commutative.
\end{theorem}
The proof strategy, carried out in Section \ref{sec:comm}, is to reduce to the {\em projective local model}: 
\begin{equation}\label{eqn:proj_local_model}
\mathbb{P}_Z(N_{Z/X}\oplus \mathcal{O}_Z)\simeq \mathbb{P}_{\mathbb{P}^r}(\mathcal{O}_{\mathbb{P}^r}^{r+1}\oplus \mathcal{O}_{\mathbb{P}^r})\dashrightarrow  \mathbb{P}(\mathcal{O}_{\mathbb{P}^r}^{r+1}\oplus \mathcal{O}_{\mathbb{P}^r})\simeq\mathbb{P}_{Z'}(N_{Z'/X'}\oplus \mathcal{O}_{Z'}). 
\end{equation}
More details about projective local models can be found in \cite{LLW}. Since the projective local model is a toric crepant wall-crossing, we can apply results of \cite{CIJ}. The reduction to the projective local model requires that we degenerate elements of $K$-groups and their Fourier-Mukai images. We achieve this in Section \ref{sec:sheaves} using results of \cite{lw}.

The simple flop $X\dashrightarrow X'$ is an example of a {\em K-equivalence} introduced in general in \cite{W}. In literature, a K-equivalence is also called a {\em crepant transformation}. It is conjectured in \cite{K} that K-equivalent varieties have equivalent derived categories. It is conjectured that varieties related by crepant transformations have equivalent Gromov-Witten theories, in a suitable sense, see \cite{CIT} and \cite{CR} for discussions in this direction. The commutativity of (\ref{diagram:descendant_FM}) can be viewed as showing that, in case of simple flops, the derived equivalence and {\em descendant} Gromov-Witten equivalence are compatible. This compatibility result has previously been proven for other crepant transformations: to our knowledge, the proven cases include crepant toric wall-crossings \cite{CIJ}, Hilbert-Chow morphism $\mathsf{Hilb}^n(\mathbb{C}^2)\to\mathsf{Sym}^n(\mathbb{C}^2)$ \cite{PHHT}, and Grassmannian flops \cite{LSW, PSW}. This paper provides one more example of such compatibility.

\subsection{Acknowledgment}
We thank Professor Chin-Lung Wang for helpful discussions and encouragement. J.-C. C. is supported by the National Science and Technology Council (NSTC) of Taiwan under Grant No. NSTC-114-2115-M-007-004-. H.-H. T. is supported in part by a Simons Foundation collaboration grant.

\section{Descendant correspondence}\label{sec:descendant}
In this section, we construct the map $$\mathbb{U}:\widetilde{\mathcal{H}}^X\to \widetilde{\mathcal{H}}^{X'}.$$ 
We work in Givental's symplectic space formalism for genus $0$ Gromov-Witten theory, see \cite{G} for details. 

\subsection{Genus zero}
By \cite{LLW}, the genus $0$ ancestor Gromov-Witten theories of $X$ and $X'$ match after analytic continuation. More precisely, generating functions of Gromov-Witten theory of $X$ (respectively $X'$) are formal power series in quantum variables\footnote{Also known as Novikov variables.} of $X$ (respectively $X'$), which index effective curve classes of $X$ (respectively $X'$). Among these quantum variables, the variable $q$ (respectively $q'$) that indexes the class\footnote{This is the class of the so-called {\em flopping curves} or {\em extremal curves} of the flop $X\dashrightarrow X'$.} of a line in $Z\subset X$ (respectively $Z'\subset X'$) plays a special role. By the results of \cite{LLW}, generating functions of Gromov-Witten theory of $X$ (respectively $X'$) are analytic near $q=0$ (respectively $q'=0$). Generating functions of Gromov-Witten theory of $X$ can be analytically continued to near $q=\infty$ and then identified with generating functions of Gromov-Witten theory of $X'$ using the additive isomorphism 
\begin{equation}\label{eqn:additive_isom}
H^*(X)\simeq H^*(X')
\end{equation}
induced by the closure of the graph of $X\dashrightarrow X'$ and suitable identifications of quantum variables, in particular $q'$ is identified with $1/q$.

In terms of Lagrangian cones, we can write the results of \cite{LLW} as
$$\mathcal{L}^{an}_{X,\tau}=\mathcal{L}^{an}_{X',\tau'},$$
after analytic continuation. Here the parameters $\tau\in H^*(X)$ and $\tau'\in H^*(X')$ are identified using (\ref{eqn:additive_isom}).

By ancestor/descendant correspondence \cite[Appendix 2]{CG}, we have the following relations between Lagrangian cones of ancestor and descendant theories:
$$\mathcal{L}_X=S_{X,\tau}^{-1}\mathcal{L}^{an}_{X,\tau} \quad\text{ and } \quad \mathcal{L}_{X'}=S_{X',\tau'}^{-1}\mathcal{L}^{an}_{X',\tau'}.$$ 
Here $S_{X,\tau}$ and $S_{X',\tau'}$ are fundamental solutions to quantum differential equations of $X$ and $X'$ given explicitly in terms of genus $0$ descendant Gromov-Witten invariants, see e.g \cite[Section 5]{G_symp}. 

The results of \cite{LLW} imply that the quantum differential equations of $X$ and $X'$ match after the analytic continuation mentioned above. By theory of linear differential equations, there exists a unique $\mathbb{U}$ that measures the difference between the fundamental solutions $S_{X,\tau}$ analytically continued to $q=\infty$ and the fundamental solution $S_{X',\tau'}$ of the quantum differential equation, i.e. 
\begin{equation}\label{eqn:U}
\mathbb{U}S_{X,\tau}^{-1}=S_{X',\tau'}^{-1}.
\end{equation}
Here, to identify $\widetilde{\mathcal{H}}^X$ and $\widetilde{\mathcal{H}}^{X'}$, we use the additive isomorphism $H^*(X)\simeq H^*(X')$ in (\ref{eqn:additive_isom}). Therefore, we obtain, after analytic continuation,
$$\mathbb{U}\mathcal{L}_X=\mathcal{L}_{X'},$$ 
which is viewed as a correspondence between genus $0$ descendant Gromov-Witten theories of $X$ and $X'$.

Furthermore, the construction implies that $\mathbb{U}$ is independent of all (quantum) variables, and we can study $\mathbb{U}$ by restricting to the one-variable case, namely keeping only the variable $q$ (and $q'$) and setting all other variables to $0$.

\subsection{Higher genus}
In \cite{ILLW}, the all-genera ancestor Gromov-Witten theories of $X$ and $X'$ are proven to match after analytic continuation. In terms of the total ancestor potentials, we can write the results of \cite{ILLW} as
\begin{equation}\label{eqn:ancestor_equal}
\mathcal{A}_{\tau}^X=\mathcal{A}_{\tau'}^{X'},    
\end{equation}
after analytic continuation. The all-genera descendant/ancestor correspondence \cite[Appendix 2]{CG} gives a relation between total descendant and ancestor potentials, written using Givental's quantization formalism:
\begin{equation}\label{eqn:anc_des_all_gen}
\mathcal{D}_X=C_X\widehat{S}_{X,\tau}^{-1}\mathcal{A}_{X,\tau}, \quad \mathcal{D}_{X'}=C_{X'}\widehat{S}_{X',\tau'}^{-1}\mathcal{A}_{X',\tau'}.
\end{equation}
Here $C_X$ and $C_{X'}$ are certain nonzero scalars explicitly expressed in terms of genus $1$ Gromov-Witten invariants. By (\ref{eqn:U}), (\ref{eqn:ancestor_equal}), and (\ref{eqn:anc_des_all_gen}), we obtain the following all-genera descendant correspondence:
\begin{equation}\label{eqn:des_equal_all_gen}
C_X^{-1}\mathcal{D}_X=\widehat{\mathbb{U}}(C_{X'}^{-1}\mathcal{D}_{X'}), 
\end{equation}
after analytic continuation.

\begin{remark}
The construction of $\mathbb{U}$ is parallel to the case of $\mathsf{Hilb/Sym}$ correspondence for $\mathbb{C}^2$ studied in \cite{PHHT}.
\end{remark}

%We try to use the result of \cite{c} in connection with degeneration formula.
\section{Deformation to the normal cone}
In this section, we discuss several constructions involving deformations to the normal cones. These constructions are used in the proof of Theorem \ref{thm:comm}.

\subsection{Specialization map of cohomology}\label{sec:specialization_map}
Consider the deformation $\mathfrak{X}'\to\mathbb{P}^1$ to the normal cone of $Z'\subset X'$. For $t\in \mathbb{P}^1\setminus\{0\}$, the fiber $\mathfrak{X}'_t$ is isomorphic to $X'$. For $t=0$, the fiber $\mathfrak{X}'_0$ is a union of $\text{Bl}_{Z'}X'$ and $\mathbb{P}_{Z'}(N_{Z'/X'}\oplus \mathcal{O}_{Z'})$ glued along $$\mathbb{P}_{Z'}(N_{Z'/X'})\subset \mathbb{P}_{Z'}(N_{Z'/X'}\oplus \mathcal{O}_{Z'}),$$ 
where $\mathbb{P}_{Z'}(N_{Z'/X'})$ is embedded into $\text{Bl}_{Z'}X'$ as the exceptional divisor. The specialization map of Borel-Moore homology groups associated to this deformation, constructed by pushing forward into a tubular neighborhood of the special fiber, followed by a retract, is {\em injective}:
\begin{equation}
H^{BM}_*(X')=H^{BM}_*(\mathfrak{X}'_t)\xhookrightarrow{\text{sp}_{\mathfrak{X}'/\mathbb{P}^1}} H^{BM}_*(\mathfrak{X}'_0).    
\end{equation}
To see this\footnote{This argument is certainly not new. We include it here for completeness.}, consider the blowup map $b:\mathfrak{X}'=\text{Bl}_{Z'\time 0}(X'\times  \mathbb{P}^1)\to X'\times \mathbb{P}^1$. $b$ is proper. For $t\in\mathbb{P}^1$, denote by $b_t:\mathfrak{X}'_t\to X'$ the restriction of $b$ to $t$. Since specialization maps are functorial with respect to proper pushforwards, we have the following commutative diagram
\begin{equation}
\xymatrix@C=5pc{
H^{BM}_*(\mathfrak{X}'_t)\ar[d]_{(b_t)_*}\ar[r]^{\text{sp}_{\mathfrak{X}'/\mathbb{P}^1}}& H^{BM}_*(\mathfrak{X}'_0)\ar[d]^{(b_0)_*}\\
H^{BM}_*({X}')\ar[r]^{\text{sp}_{{X}'\times \mathbb{P}^1/\mathbb{P}^1}}& H^{BM}_*({X}').
}    
\end{equation}
Since $(b_t)_*$ and $\text{sp}_{{X}'\times \mathbb{P}^1/\mathbb{P}^1}$ are identities, the injectivity of $\text{sp}_{\mathfrak{X}'/\mathbb{P}^1}$ follows.

\iffalse
In fact, the construction shows that injectivity holds for the specialization map associated to $\tilde{\mathfrak{X}'}$:
\begin{equation}
H^{BM}_*(X')=H^{BM}_*(\tilde{\mathfrak{X}'}_t)\hookrightarrow H^{BM}_*(\tilde{\mathfrak{X}'}_0).    
\end{equation}
\fi
In cohomology, we have 
\begin{equation}
H^*(X')\underset{\text{duality}}{\simeq} H^{BM}_*(X')\hookrightarrow H^{BM}_*({\mathfrak{X}'}_0)\underset{\text{duality}}{\simeq} H^*({\mathfrak{X}'}_0)\to H^*(\text{Bl}_{Z'}X')\oplus H^*(\mathbb{P}_{Z'}(N_{Z'/X'}\oplus \mathcal{O}_{Z'})),
\end{equation}
where the last map is part of the Mayer-Vietoris sequence. The kernel of the last map is the image of $H^{*-1}(\mathbb{P}_{Z'}(N_{Z'/X'}))\to H^*({\mathfrak{X}'}_0)$. As we work with even degree cohomology\footnote{The restriction to even degree cohomology is reasonable, because the ingredients involved in Theorem \ref{thm:comm} are all algebraic.}, this kernel is $0$. This gives an embedding
\begin{equation}
\Phi': H^*(X')\to H^*(\text{Bl}_{Z'}X')\oplus H^*(\mathbb{P}_{Z'}(N_{Z'/X'}\oplus \mathcal{O}_{Z'}))
\end{equation}
and the analogous $\Phi$ for $H^*(X)$.

\subsection{Degenerations of sheaves}\label{sec:sheaves}
 Denote by $$\pi:\mathfrak{X}\to \mathbb{P}^1$$ the deformation to the normal cone of $Z\subset X$. For $t\in \mathbb{P}^1\setminus\{0\}$, the fiber $\mathfrak{X}_t$ is isomorphic to $X$. For $t=0$, the fiber $\mathfrak{X}_0$ is a union of $\text{Bl}_{Z}X$ and $\mathbb{P}_Z(N_{Z/X}\oplus \mathcal{O}_Z)$ glued along $\mathbb{P}_Z(N_{Z/X})\subset \mathbb{P}_Z(N_{Z/X}\oplus \mathcal{O}_Z)$ embedded into $\text{Bl}_ZX$ as the exceptional divisor. There is a divisor $$\mathfrak{Z}\subset \mathfrak{X}/\mathbb{P}^1$$ that restricts to $Z$ and to the other section of $\mathbb{P}_Z(N_{Z/X}\oplus \mathcal{O}_Z)$.

Pick a $\pi$-ample line bundle $\mathcal{L}$ on $\mathfrak{X}$. 

Let $\mathcal{F}$ be a coherent sheaf on $X$. Pick $M, N\gg 1$ such that there is a surjection $\mathcal{O}_X^{\oplus M}\to \mathcal{F}\otimes \mathcal{L}_t^{N}$. So 
\begin{equation}
(\mathcal{L}_t^{-N})^{\oplus M}\to \mathcal{F}
\end{equation}
is a point in a Quot scheme of $X$. Let $P$ be the Hilbert polynomial of this point. 

Consider the Quot scheme $$\mathfrak{Quot}^{(\mathcal{L}^{-N})^{\oplus M}, P}(\mathfrak{X}/\mathbb{P}^1).$$ 
The quotient $(\mathcal{L}^{-N})^{\oplus M}|_{\mathbb{P}^1\setminus\{0\}}\to \mathcal{F}\boxtimes \mathcal{O}_{\mathbb{P}^1\setminus\{0\}}$ gives a non-constant map $$\mathbb{P}^1\setminus\{0\}\to \mathfrak{Quot}^{(\mathcal{L}^{-N})^{\oplus M}, P}(\mathfrak{X}/\mathbb{P}^1).$$ By the properness result of \cite[Sections 4 and 5]{lw}, this extends to a morphism $\mathbb{P}^1\to \mathfrak{Quot}^{(\mathcal{L}^{-N})^{\oplus M}, P}(\mathfrak{X}/\mathbb{P}^1)$, which is a flat family of quotients on a {\em modification}\footnote{This modification depends on the sheaf $\mathcal{F}$.} $$\tilde{\mathfrak{X}}\to\mathbb{P}^1$$ of $\mathfrak{X}\to\mathbb{P}^1$. This gives a quotient  $$(\mathcal{L}^{-N})^{\oplus M}\to\mathfrak{F}$$
on $\tilde{\mathfrak{X}}$. 
On the fiber over $0\in\mathbb{P}^1$, this gives a quotient 
\begin{equation}
v^*(\mathcal{L}_0^{-N})^{\oplus M}\to \mathcal{F}_0.    
\end{equation}
Here $v:\mathfrak{X}_0[n]\to\mathfrak{X}_0$ is a suitable expansion. We know that $$\mathfrak{X}_0[n]=\text{Bl}_ZX\coprod C[n]\coprod \mathbb{P}_Z(N_{Z/X}\oplus \mathcal{O}_Z)$$ glued along suitable divisors. Here $C[n]$, the ``rubber'', is a chain of $\mathbb{P}^1$-bundles, see \cite{lw} for details. The coherent sheaf $\mathcal{F}_0$ is obtained from coherent sheaves on $\text{Bl}_ZX$, $C[n]$, and $\mathbb{P}_Z(N_{Z/X}\oplus \mathcal{O}_Z)$ with compatibility conditions on gluing divisors. By construction, the sheaf on $C[n]$ is invariant with respect to the $\mathbb{G}_m^n$ action on $C[n]$. Thus the sheaf on $C[n]$ is determined by those on $\text{Bl}_ZX$ and $\mathbb{P}_Z(N_{Z/X}\oplus \mathcal{O}_Z)$.

The divisor $\mathfrak{Z}$ lifts to a divisor $\tilde{\mathfrak{Z}}\subset \tilde{\mathfrak{X}}$. The blow-up 
$$\pi_{\tilde{\mathfrak{X}}}:\text{Bl}_{\tilde{\mathfrak{Z}}} \tilde{\mathfrak{X}}\to \tilde{\mathfrak{X}}$$
is a $\mathbb{P}^1$-family with exceptional divisor a $\mathbb{P}^1$-family whose fibers are isomorphic to $\mathbb{P}^r\times\mathbb{P}^r$. Blowing down the other ruling yields a $\mathbb{P}^1$-family 
$$\pi_{\tilde{\mathfrak{X}}'}:\text{Bl}_{\tilde{\mathfrak{Z}}} \tilde{\mathfrak{X}}\to\tilde{\mathfrak{X}'}.$$
Over $t\neq 0$, the fiber $\tilde{\mathfrak{X}'}_t$ is isomorphic to $X'$. The fiber over $0$ is $$\mathfrak{X}'_0[n]=\text{Bl}_{Z'}X'\coprod C[n]\coprod \mathbb{P}_{Z'}(N_{Z'/X'}\oplus \mathcal{O}_{Z'})$$ 
glued along suitable divisors. Furthermore, the birational map 
$$\tilde{\mathfrak{X}}\dashrightarrow \tilde{\mathfrak{X}}'$$ 
is a flat $\mathbb{P}^1$-family that deforms the flop $X\dashrightarrow X'$ to the birational map 
\begin{equation}\label{eqn:flop_0}
\text{Bl}_ZX\coprod C[n]\coprod \mathbb{P}_Z(N_{Z/X}\oplus \mathcal{O}_Z)\dashrightarrow \text{Bl}_{Z'}X'\coprod C[n]\coprod \mathbb{P}_{Z'}(N_{Z'/X'}\oplus \mathcal{O}_{Z'})
\end{equation}
which is the identity on the first two factors and the projective local model of the simple flop $X\dashrightarrow X'$ on the third factor.

\section{Commutativity: proof of the main result}\label{sec:comm}
In this section, we prove Theorem \ref{thm:comm}.

For $\alpha\in K(X)$ represented by a coherent sheaf $\mathcal{F}$, we want to show that 
\begin{equation}\label{eqn:main_question}
\Psi^{X'}(\mathbb{FM}(\alpha))=\mathbb{U}\Psi^X(\alpha). 
\end{equation}
This is an equality in $\widetilde{\mathcal{H}}^{X'}=H^*(X')[[\log z]]((z^{-1}))$.  

By assumption, $\mathbb{FM}(\alpha)$ can be represented by a complex of coherent sheaves on $X'$. By definition of $\Psi^X$ and $\Psi^{X'}$, $\Psi^X(\alpha)$ and $\Psi^{X'}(\mathbb{FM}(\alpha))$ take values $\text{CH}^*(X)[[\log z]]((z^{-1}))$ and $\text{CH}^*(X')[[\log z]]((z^{-1}))$ respectively. Namely they have coefficients in Chow groups. Hence, their coefficents can be made to take values in Borel-Moore homology by cycle maps, and then in (even degree) cohomology by duality.

By the injectivity discussed in Section \ref{sec:specialization_map}, it suffices to embed the equality (\ref{eqn:main_question}) into $$(H^*(\text{Bl}_{Z'}X')\oplus H^*(\mathbb{P}_{Z'}(N_{Z'/X'}\oplus \mathcal{O}_{Z'})))[[\log z]]((z^{-1}))$$ 
and prove it there.

Since the specialization map for Chow groups is obtained by taking closures of cycles, the existence of the sheaf $\mathfrak{F}$ on $\tilde{\mathfrak{X}}$ associated to $\mathcal{F}$, as constructed in Section \ref{sec:sheaves}, implies that the image of $\text{ch}(\mathfrak{F}|_{\tilde{\mathfrak{X}}_t})$ under the specialization map is the intersection of the closure of $\text{ch}(\mathfrak{F}|_{\tilde{\mathfrak{X}}\setminus \tilde{\mathfrak{X}}_0})$ with $\tilde{\mathfrak{X}}_0$. In other words, it is $\text{ch}(\mathfrak{F}|_{\tilde{\mathfrak{X}}_0})$. Hence, 
\begin{equation}
\text{ch}(\mathcal{F}_0|_{\text{Bl}_ZX})+\text{ch}(\mathcal{F}_0|_{\mathbb{P}_Z(N_{Z/X}\oplus \mathcal{O}_Z)})=\Phi(\text{ch}(\mathcal{F})),    
\end{equation}

By construction, the specialization of the cycle $\text{ch}(\mathbb{FM}(\mathcal{F}))$ is represented by the closure of the class $\text{ch}(R(\pi_{\mathfrak{X}'})_*L\pi_{\mathfrak{X}}^*(\mathcal{F}\boxtimes \mathcal{O}_{\mathbb{P}^1\setminus\{0\}}))$ in $\tilde{\mathfrak{X}}'_0$. Therefore, we have
\begin{equation}
\text{ch}(\mathbb{FM}(\mathcal{F}_0)|_{\text{Bl}_{Z'}X'})+\text{ch}(\mathbb{FM}(\mathcal{F}_0)|_{\mathbb{P}_{Z'}(N_{Z'/X'}\oplus \mathcal{O}_{Z'})})=\Phi'(\text{ch}(\mathbb{FM}(\mathcal{F}))).    
\end{equation}
Since the birational map (\ref{eqn:flop_0}) is the identity on the first two factors, we have $$\text{ch}(\mathbb{FM}(\mathcal{F}_0)|_{\text{Bl}_{Z'}X'})=\text{ch}(\mathcal{F}_0)|_{\text{Bl}_{Z}X}.$$

Applying the above argument to the sheaf $\Omega_{\tilde{\mathfrak{X}}/\mathbb{P}^1}^\vee$, we have 
\begin{equation}
\Phi(\Gamma(T_X))=\Gamma(T_{\text{Bl}_ZX})+\Gamma(T_{\mathbb{P}_Z(N_{Z/X}\oplus \mathcal{O}_Z)}),     \end{equation}
and similarly for $\Gamma(T_{X'})$. Also,
\begin{equation}
\Phi(c_1(T_X))=c_1(T_{\text{Bl}_ZX})+c_1(T_{\mathbb{P}_Z(N_{Z/X}\oplus \mathcal{O}_Z)}),
\end{equation}
and similarly for $c_1(T_{X'})$. This implies the compatibility of $\rho^X$ (respectively  $\rho^{X'}$) with $\Phi$ (respectively $\Phi'$). 

The compatibility of $\mu^X$ and $\text{deg}_0^X$ (respectively $\mu^{X'}$ and $\text{deg}_0^{X'}$) with $\Phi$ (respectively $\Phi'$) concerns grading and is straightforward:
\begin{equation}
\Phi\circ \mu^X=(\mu^{\text{Bl}_ZX}\oplus \mu^{\mathbb{P}_Z(N_{Z/X}\oplus \mathcal{O}_Z)})\circ\Phi,\quad \Phi\circ \text{deg}_0^X= (\text{deg}_0^{\text{Bl}_ZX}\oplus \text{deg}_0^{\mathbb{P}_Z(N_{Z/X}\oplus \mathcal{O}_Z)})\circ\Phi,\end{equation}
and similarly for $X'$.

The map $\mathbb{U}$ is constructed to satisfy 
\begin{equation}\label{eqn:fund_soln}
\mathbb{U}S_{X,\tau}^{-1}=S_{X',\tau'}^{-1}.
\end{equation}
We need to study (\ref{eqn:fund_soln}) under the embeddings $\Phi$ and $\Phi'$.

As noted above, $\mathbb{U}$ does not depend on any quantum variable. Therefore, in (\ref{eqn:fund_soln}), we can set to $0$ all quantum variables other than $q$ and $q'$ (which index the flopping curves). Then, Gromov-Witten invariants of $X$ appearing in (\ref{eqn:fund_soln}) are of the form $$\langle \gamma_1,...,\gamma_k/(z-\psi)\rangle_{0,\text{extremal}}^X,$$ namely, those Gromov-Witten invariants whose degrees are multiples of the class of a line in $Z\subset X$. The flopping curves in $X$ are contained in $Z\subset X$, with normal bundle $N_{Z/X}\simeq \mathcal{O}_{\mathbb{P}^r}(-1)^{\oplus r+1}$. Therefore, flopping curves are constrained to lie in the divisor $\tilde{\mathfrak{Z}}$ in the degeneration $\tilde{\mathfrak{X}}\to \mathbb{P}^1$. The divisor $\tilde{\mathfrak{Z}}$ meets the special fiber $\tilde{\mathfrak{X}}_0$ only at the factor $\mathbb{P}_{Z}(N_{Z/X}\oplus \mathcal{O}_{Z})$ and is disjoint from the section $\mathbb{P}_Z(N_{Z/X})$ where the two pieces $\text{Bl}_{Z}X$ and $\mathbb{P}_Z(N_{Z/X}\oplus \mathcal{O}_Z)$ glued. This implies that only the factor $\mathbb{P}_Z(N_{Z/X}\oplus \mathcal{O}_Z)$ contributes when evaluating $\langle \gamma_1,...,\gamma_k/(z-\psi)\rangle_{0,\text{extremal}}^X$ by applying the degeneration formula \cite{Li} associated to the family $\mathfrak{X}\to \mathbb{P}^1$. We can also choose $\tau\in H^*(X)$ so that the $H^*(\mathbb{P}_Z(N_{Z/X}\oplus \mathcal{O}_Z))$-summand of $\Phi(\tau)$ is not zero. Then, when restricted to the $q$-direction, applying the degeneration formula for Gromov-Witten invariants to $\mathfrak{X}\to \mathbb{P}^1$ implies that the solution $S_{X,\tau}^{-1}$ satisfies
 \begin{equation*}
\Phi S_{X,\tau}^{-1}|_{\text{extremal}}=\begin{pmatrix}
\text{Id}_{\text{Bl}_ZX}& 0\\
0 & S_{\mathbb{P},\tau_{\mathbb{P}}}^{-1}|_{\text{extremal}} \end{pmatrix}\cdot \Phi.
 \end{equation*}
Here, $|_{\text{extremal}}$ indicates setting all quantum variables to $0$ except $q$. $$\mathbb{P}_{Z}(N_{Z/X}\oplus \mathcal{O}_{Z})=\mathbb{P}\dashrightarrow \mathbb{P}'=\mathbb{P}_{Z'}(N_{Z'/X'}\oplus \mathcal{O}_{Z'})$$
is the projective local model described in (\ref{eqn:proj_local_model}). $S_{\mathbb{P},\tau_{\mathbb{P}}}$ is the corresponding fundamental solution to the quantum differential equation for $\mathbb{P}$. $\tau_{\mathbb{P}}\in H^*(\mathbb{P})$ is the $H^*(\mathbb{P})$-summand of $\Phi(\tau)$. 
 
By a similar argument, when restricted to the $q'$-direction, we have 
\begin{equation*}
\Phi'S_{X',\tau'}^{-1}|_{\text{extremal}}=\begin{pmatrix}
\text{Id}_{\text{Bl}_{Z'}X'}& 0\\
0 & S_{\mathbb{P}',\tau_{\mathbb{P}'}}^{-1}|_{\text{extremal}} 
\end{pmatrix}\cdot \Phi'.
\end{equation*}
Here, $|_{\text{extremal}}$ indicates setting all quantum variables to $0$ except $q'$. $S_{\mathbb{P}',\tau_{\mathbb{P}'}}$ is the corresponding fundamental solution to the quantum differential equation for $\mathbb{P}'$. $\tau_{\mathbb{P}'}\in H^*(\mathbb{P}')$ is the $H^*(\mathbb{P}')$-summand of $\Phi'(\tau')$.  

It follows that, after setting all quantum variables to $0$ except $q$ and $q'$, (\ref{eqn:fund_soln}) can be written as 
 \begin{equation}
\begin{pmatrix}
\text{Id}_{\text{Bl}_ZX} &0 \\
0 & \mathbb{U}_{\mathbb{P}}|_{\text{extremal}}\end{pmatrix}     
\begin{pmatrix}
\text{Id}_{\text{Bl}_ZX} & 0\\
0 & S_{\mathbb{P},\tau_{\mathbb{P}}}^{-1}|_{\text{extremal}}
\end{pmatrix}
=
\begin{pmatrix}
\text{Id}_{\text{Bl}_{Z'}X'} & 0\\
0 & S_{\mathbb{P}',\tau_{\mathbb{P}'}}^{-1}|_{\text{extremal}}
\end{pmatrix}
\end{equation}
after applying the embeddings $\Phi, \Phi'$. Namely,
\begin{equation}
\Phi'\mathbb{U}=\begin{pmatrix}
\text{Id}_{\text{Bl}_ZX} &0 \\
0 & \mathbb{U}_{\mathbb{P}}
\end{pmatrix}\cdot \Phi. 
\end{equation}

The equality $\mathbb{U}_\mathbb{P}S_{\mathbb{P},\tau_{\mathbb{P}}}^{-1}|_{\text{extremal}}=S_{\mathbb{P}',\tau_{\mathbb{P}'}}^{-1}|_{\text{extremal}}$ means that $\mathbb{U}_\mathbb{P}$ is the symplectic transformation for the projective local model. Since the projective local model is toric, $\mathbb{U}_\mathbb{P}$ is constructed as a special case of \cite[Theorem 6.1]{CIJ}. 

Putting everything together, we have 
\begin{equation}
\begin{split}
\Phi'\Psi^{X'}(\mathbb{FM}(\alpha))
&=\Psi^{\text{Bl}_{Z'}X'}(\mathbb{FM}(\mathcal{F}_0)|_{\text{Bl}_{Z'}X'})+\Psi^{\mathbb{P}_{Z'}(N_{Z'/X'}\oplus \mathcal{O}_{Z'})}(\mathbb{FM}(\mathcal{F}_0)|_{\mathbb{P}_{Z'}(N_{Z'/X'}\oplus \mathcal{O}_{Z'})})\\
&=\Psi^{\text{Bl}_{Z'}X'}(\mathbb{FM}(\mathcal{F}_0|_{\text{Bl}_{Z}X}))+\Psi^{\mathbb{P}_{Z'}(N_{Z'/X'}\oplus \mathcal{O}_{Z'})}(\mathbb{FM}(\mathcal{F}_0|_{\mathbb{P}_{Z}(N_{Z/X}\oplus \mathcal{O}_{Z})}))\\
&=\Psi^{\text{Bl}_{Z}X}(\mathcal{F}_0|_{\text{Bl}_{Z}X})+\Psi^{\mathbb{P}_{Z'}(N_{Z'/X'}\oplus \mathcal{O}_{Z'})}(\mathbb{FM}(\mathcal{F}_0|_{\mathbb{P}_{Z}(N_{Z/X}\oplus \mathcal{O}_{Z})}))
\end{split}    
\end{equation}
and 
\begin{equation}
\begin{split}
\Phi'\mathbb{U}\Psi^X(\alpha)
&=\begin{pmatrix}
\text{Id} &0 \\
0 & \mathbb{U}_{\mathbb{P}}\end{pmatrix}\cdot \Phi\Psi^X(\alpha)\\
&=\begin{pmatrix}
\text{Id} &0 \\
0 & \mathbb{U}_{\mathbb{P}}\end{pmatrix}(\Psi^{\text{Bl}_ZX}(\mathcal{F}_0|_{\text{Bl}_ZX})+\Psi^{\mathbb{P}_Z(N_{Z/X}\oplus \mathcal{O}_Z)}(\mathcal{F}_0|_{\mathbb{P}_Z(N_{Z/X}\oplus \mathcal{O}_Z)}))\\
&=\Psi^{\text{Bl}_{Z}X}(\mathcal{F}_0|_{\text{Bl}_{Z}X})+\mathbb{U}_{\mathbb{P}}\Psi^{\mathbb{P}_Z(N_{Z/X}\oplus \mathcal{O}_Z)}(\mathcal{F}_0|_{\mathbb{P}_Z(N_{Z/X}\oplus \mathcal{O}_Z)}).
\end{split}
\end{equation}
Since $\Phi$ and $\Phi'$ are embeddings, (\ref{eqn:main_question}) follows from 
\begin{equation}
\Psi^{\mathbb{P}_{Z'}(N_{Z'/X'}\oplus \mathcal{O}_{Z'})}(\mathbb{FM}(\mathcal{F}_0|_{\mathbb{P}_{Z'}(N_{Z'/X'}\oplus \mathcal{O}_{Z'})}))=\mathbb{U}_{\mathbb{P}}\Psi^{\mathbb{P}_Z(N_{Z/X}\oplus \mathcal{O}_Z)}(\mathcal{F}_0|_{\mathbb{P}_Z(N_{Z/X}\oplus \mathcal{O}_Z)}),
\end{equation}
which is a special case of \cite[Theorem 6.1]{CIJ} concerning crepant toric wall-crossings. This completes the proof.

\begin{remark}
A class of K-equivalence generalizing simple flops is {\em ordinary flops} \cite{LLW1}. For an ordinary flop, the exceptional loci are projective {\em bundles} over a smooth base. The invariance of genus $0$ Gromov-Witten theory under ordinary flops has been proved in \cite{LLW2, LLQW}. The projective local model of an ordinary flop is a double projective bundle, which is an example of toric bundles in the sense of \cite{SU}, \cite{CGT}, \cite{Koto}. For {\em split} toric bundles, descendant correspondence between genus $0$ Gromov-Witten theory under crepant toric wall-crossings has been studied in \cite{CCT}. An extension of \cite{CCT} to non-split toric bundles should be possible, in light of \cite{Koto}. Therefore, the approach in this paper can be used to construct the diagram (\ref{diagram:descendant_FM}) for ordinary flops of split type and prove its commutativity. In this paper, we want to emphasize our approach to the diagram (\ref{diagram:descendant_FM}) and do not pursue the more general case.
\end{remark}


\begin{thebibliography}{12}

\bibitem{BO} A. Bondal, D. Orlov, {\em Semiorthogonal decomposition for algebraic varieties}, arXiv:alg-geom/9506012.

\bibitem{CCT} Q. Chao, J.-C. Chen, H.-H. Tseng,{\em Crepant transformation correspondence for toric stack bundles}, Adv. Geom. 26, No. 1, 105--125 (2026).

%\bibitem{c} T. Coates, {\em Givental’s Lagrangian cone and $S^1$–equivariant Gromov-Witten theory}, Math. Res. Lett. 15, No. 1, 15--31 (2008).

\bibitem{CG} T. Coates, A. Givental, {\em Quantum Riemann-Roch, Lefschetz and Serre}, Ann. Math. (2) 165, No. 1, 15--53 (2007).

\bibitem{CGT} T. Coates, A. Givental, H.-H. Tseng, {\em Virasoro constraints for toric bundles}, Forum Math. Pi 12, Paper No. e4, 28 p. (2024).

\bibitem{CIJ} T. Coates, H. Iritani, Y. Jiang, {\em The crepant transformation conjecture for toric complete intersections}, Adv. Math. 329, 1002--1087 (2018).

\bibitem{CIT} T. Coates, H. Iritani, H.-H. Tseng, {\em Wall-crossings in toric Gromov-Witten theory. I: Crepant examples}, Geom. Topol. 13, No. 5, 2675--2744 (2009).

\bibitem{CR} T. Coates, Y. Ruan, {\em Quantum cohomology and crepant resolutions: a conjecture. (Cohomologie quantique et résolutions crépantes : une conjecture.)}, Ann. Inst. Fourier 63, No. 2, 431--478 (2013).

%\bibitem{fulton} W. Fulton, {\em Intersection theory},

\bibitem{G_symp} A. Givental, {\em Gromov-Witten invariants and quantization of quadratic Hamiltonians}, Mosc. Math. J.  {\bf 4} (2001), 551--568.

\bibitem{G} A. Givental, {\em Symplectic geometry of Frobenius structures}, in: ``Frobenius manifolds. Quantum cohomology and singularities'', 91--112, Aspects of Mathematics E 36, Vieweg (2004).

\bibitem{I} H. Iritani, {\em An integral structure in quantum cohomology and mirror symmetry for toric orbifolds}, Adv. Math. 222, No. 3, 1016--1079 (2009).

\bibitem{ILLW} Y. Iwao, Y.-P. Lee, H.-W. Lin, C.-L. Wang, {\em Invariance of Gromov-Witten theory under a simple flop}, J. Reine Angew. Math. 663, 67--90 (2012).

\bibitem{K} Y. Kawamata, {\em D-equivalence and K-equivalence}, J. Differ. Geom. 61, No. 1, 147--171 (2002).

\bibitem{Koto} Y. Koto, {\em A mirror theorem for non-split toric bundles}, Math. Ann. 393, No. 3-4, 3337--3394 (2025).

\bibitem{LLQW} Y.-P. Lee, H.-W. Lin, F. Qu, C.-L. Wang, {\em Invariance of quantum rings under ordinary flops. III: A quantum splitting principle}, Camb. J. Math. 4, No. 3, 333--401 (2016).

\bibitem{LLW} Y.-P. Lee, H.-W. Lin, C.-L. Wang, {\em Flops, motives, and invariance of quantum rings}, Ann. Math. (2) 172, No. 1, 243--290 (2010).

\bibitem{LLW1} Y.-P. Lee, H.-W. Lin, C.-L. Wang, {\em Invariance of quantum rings under ordinary flops. I: Quantum corrections and reduction to local models}, Algebr. Geom. 3, No. 5, 578--614 (2016).

\bibitem{LLW2} Y.-P. Lee, H.-W. Lin, C.-L. Wang, {\em Invariance of quantum rings under ordinary flops. II: A quantum Leray-Hirsch theorem}, Algebr. Geom. 3, No. 5, 615--653 (2016).

\bibitem{Li} J. Li, {\em A degeneration formula of GW-invariants}, J. Differ. Geom. 60, No. 2, 199--293 (2002).

\bibitem{lw} J. Li, B. Wu, {\em Good degeneration of Quot-schemes and coherent systems}, Commun. Anal. Geom. 23, No. 4, 841--921 (2015).

\bibitem{LSW} W. Lutz, Q. Shafi, R. Webb, {\em Crepant Transformation Conjecture For the Grassmannian Flop}, Trans. Am. Math. Soc. 378, No. 7, 4671--4706 (2025), arXiv:2404.12302.

\bibitem{PHHT} R. Pandharipande, H.-H. Tseng, {\em The $\mathsf{Hilb}/\mathsf{Sym}$ correspondence
for $\mathbb{C}^2$: descendents and Fourier-Mukai}, Math. Ann. {\bf 375} (2019), 509--540.

\bibitem{PSW} N. Priddis, M. Shoemaker, Y. Wen, {\em Wall crossing and the Fourier-Mukai transform for Grassmann flops}, SIGMA 21 (2025), 008, 33 pages, arXiv:2404.12303. 

\bibitem{SU} P. Sankaran, V. Uma, {\em Cohomology of toric bundles}, Comment. Math. Helv. 78, No. 3, 540--554 (2003); errata ibid. 79, 840--841 (2004).

\bibitem{W} C.-L. Wang, {\em K-equivalence in birational geometry}, in: ``Second international congress of Chinese mathematicians. Proceedings of the congress (ICCM2001), Taipei, Taiwan, December 17–22, 2001'', New Studies in Advanced Mathematics 4, 199--216, International Press (2004).
    
\end{thebibliography}
\end{document}